\documentclass[letterpaper,12pt,leqno,oneside]{article}
\usepackage{color}
\usepackage{bm}
\usepackage{exscale}
\usepackage{amsmath}
\usepackage{amsfonts}
\usepackage{pagecolor,lipsum}
\usepackage{xcolor}
\usepackage{textcomp}
\usepackage{wasysym}
\usepackage{stmaryrd}
\usepackage{amscd}
\usepackage{graphicx}
\usepackage{amsxtra}
\usepackage{amssymb}
\usepackage{theorem}
\usepackage[final]{epsfig}
\usepackage{eqnarray}
\newtheorem{proposition}{Proposition}[subsection]

\newtheorem{lemma}[proposition]{Lemma}
{\theorembodyfont{\rmfamily}\newtheorem{remark}[proposition]{Remark}}

\newtheorem{corollary}[proposition]{Corollary}

{\theorembodyfont{\rmfamily}}

\newfont{\abc}{cmtt10 scaled 1200}

\def\R{\mathbb{R}}

\def\Z{\mathbb{Z}}

\def\ve{\varepsilon}

\def\ra{\rightarrow}
\def\cs{\symbol{35}}
\def\p{\partial}
\def\qed{\hfill $\Box$ \\}
\def\mm{\mbox}
\def\v{= \emptyset}
\def\n{\neq \emptyset}
\def\D{\mathbf{ID}}

\def\bp{\langle A \rangle}

\def\db{d^{\,\flat}}

\def\llb{\llbracket}
\def\rrb{\rrbracket}

\def\bp{\langle A \rangle}
\def\bsp{[A]}

\begin{document}
 \vspace*{0cm}

\begin{center}\Large{\bf{The Higher Dimensional Positive Mass Theorem I}}\\
\smallskip
{\small{by}}\\
\smallskip
\large{\bf{Joachim Lohkamp}}
\end{center}

\vspace{0.6cm}
\noindent Mathematisches Institut, Universit\"at M\"unster, Einsteinstrasse 62, Germany\\
 {\small{\emph{e-mail:  j.lohkamp@uni-muenster.de}}}
\vspace{0.6cm}

 {\small{\center \tableofcontents}

{\contentsline {subsection}{\numberline {}References}{. .}}}

\bigskip
%{\footnotesize \center \tableofcontents}

\setcounter{section}{1}
\renewcommand{\thesubsection}{\thesection}
\subsection{Introduction} \label{introduction}
\medskip

The positive mass conjecture asserts that an isolated relativistic gravitational system, like a galaxy, has a positive total mass, unless it is entirely vacuous. Then, and only in this case,  the mass vanishes. \\

The conjecture supports the physical plausibility of general relativity since, otherwise, there were finite gravitating systems able to radiate arbitrarily large amounts of energy.  This would cause an instability of gravitational systems. Mathematically, this is a delicate conjecture, in particular, because the equivalence principle shows that there is no local measure for the energy content of the gravitational field. \\

Since most of the attempts to unify all four basic physical interactions assume the existence of extra dimensions, the positive mass conjecture is a valuable challenge to check their  physical plausibility and remains relevant in  dimensions $\ge 5$.\\

A key case of the conjecture is that of time-symmetric space-times. This is the so-called Riemannian positive mass conjecture. It is a model case for more general versions and it plays an important role in scalar curvature geometry, again, in arbitrary dimensions.\\

Concretely, it is equivalent to the assertion, explicitly formulated  in Theorem 2  below, that there is no general mechanism to locally deform a manifold in such a way that it increases the scalar curvature, even when we admit topological modifications. This underlines a basic asymmetry: we can locally and arbitrarily strongly \emph{decrease} the scalar curvature from arbitrarily small changes of the metric geometry, [L1],[L2].\\

Also, the positivity of mass is an essential ingredient in the study of various aspects of the Yamabe problem claiming that any metric on a compact manifold is conformal to one with constant scalar curvature. The impact reaches deeply from the existence to compactness theorems for solutions of this problem, [LP],[KMS].\\

\textbf{Classical and New Results} \, The  Riemannian version of the positive mass conjecture, like its space-time siblings, for asymptotically flat spaces, has been established in the cases where either the dimension of the space is $\le 7$, by Schoen and Yau using minimal hypersurface techniques [SY1]-[SY3],[S],  and recent accounts [E] and [EHLS],  or the underlying manifold carries a spin topology. This latter approach is due to Witten and it uses Bochner type arguments for spinors [W], [B], [PT].\\

 In this paper we prove the Riemannian positive mass conjecture in all dimensions and without any topological constraints.\\

Our approach generalizes the minimal hypersurface approach, introduced by Schoen and Yau, which properly worked only for regular hypersurfaces. This, in turn, was the reason why their argument applied only in low dimensions, cf. Ch.\ref{int2} below.\\

To treat the general case, we  incorporate also singular hypersurfaces. Then, we apply dedicated surgery techniques,  we developed in [L3], to eliminate the  singularities from the picture. This regularization scheme essentially employs the hyperbolic unfoldings of such hypersurfaces and, from this, the fine control over the elliptic analysis on such hypersurfaces towards their singularities,  we established in [L1] and [L2], in terms of so-called skin structures.\\

A broader overview of the classical issues, the new methods and of the arguments in this paper, is given in Ch.\ref{int2} and Ch.\ref{sk1} below.\\\\

\subsubsection{Main Results of the Paper} \label{int1}
\bigskip

In this paper we exclusively consider the Riemannian versions of gravitational systems and their total mass.  The most basic mathematical model of an isolated gravitational system is that of an asymptotically flat geometry, [B] and [LP]. \\

\textbf{Definition 1} \textbf{(Asymptotic Flat Spaces)}\label{af} \, \emph{A complete Riemannian manifold $(M^n, g)$  is called \textbf{asymptotically flat of order} $\tau$, for some $\tau >0$, if there exists a decomposition $M = M_0\cup M_\infty$, with $M_0$ \emph{compact}, and $M_\infty$ called the \textbf{end} of $M$, so that:}\\

\emph{The end $(M_\infty, g)$ is isometric to $(\R^n \setminus B_R(0),g^*)$ for some metric $g^*$ with
\[ g^*_{ij} = \delta_{ij} + O (|x|^{- \tau}), \; \frac{\p g^*_{ij}}{\p x_k} = O (|x|^{- \tau - 1}),\; \frac{\p^2 g^*_{ij}}{\p x_k \p x_l} = O(|x|^{- \tau - 2}), \] for some $R>0$.}\\

If the asymptotically flat end $M_{\infty}$ has order $\tau > \frac{n-2}{2}$ we can assign it an invariant $E (M,g)$, the \emph{total energy} or \emph{ADM energy}, named after Arnowitt, Deser and Misner who brought the energy definition into this generally accepted form, [ADM].\\

$E (M,g)$ measures the asymptotic behaviour of $(M, g)$ near infinity. Intuitively, $E < 0$ respectively $E>0$ mean that $M$ has
some hyperbolic respectively parabolic flavour near infinity: the volume of large distance balls asymptotically grows slightly stronger respectively weaker than that of balls in the Euclidean space.\\

\textbf{Definition 2} \textbf{(Total Energy)}\label{en} \, \emph{Given an asymptotically flat end $M_{\infty}$ of order $\tau > \frac{n-2}{2}$ we  assign it the \textbf{total energy} $E (M,g)$
\begin{equation}\label{adm}
E (M,g) = \frac{1}{Vol (S^{n-1})} \cdot \lim_{R \ra
\infty} \int_{\p B_R} \sum_{i,j} \left(\frac{\p g_{ij}}{\p x_i} - \frac{\p g_{ij}}{\p x_j} \right) \cdot \nu_j \: dV_{n-1},
\end{equation}
where $\nu = (\nu_1\ldots \nu_n)$ is the outer normal vector to $\p B_R$. }\\

\textbf{Remark 1 (Mass, Energy and Scalar Curvature)} \, The expression in (\ref{adm}) is independent of the chosen asymptotic coordinates $x_i$, [B]. It can be motivated directly from Einstein's field equations applying divergence theorem arguments.\\

In subsequent accounts,cf.[L6], we also extend the proof of the Riemannian positive mass conjecture to its more general physical counterparts for space-times where one makes a finer distinction between positive \emph{energy} and  positive \emph{mass} theorems. The assertion $E \ge 0$ is actually commonly called the positive energy theorem, whereas, the sharper version $E \ge  |P|$, where $|P|$ is the total momentum, is named the positive mass theorem.\\

In the Riemannian versions, one may interpret as the time-symmetric case of the space-time version, we have $P=0$. Thus, the Riemannian form of the positive energy and the positive mass conjecture coincide.\\

 Also, a basic assumption in the physical setup for the positive mass conjectures  is the validity of the so-called \emph{dominant energy condition} (DEC). In simple terms, the DEC imposes the plausible constraint on the stress-energy tensor in Einstein's field equations that the energy must not flow faster than light and the local mass-energy density must not be negative. In the Riemannian case, this translates into the assumption that the \emph{scalar curvature} $scal (g)$\emph{ is non-negative}, $scal (g) \ge 0$, [LP],Ch.8.\qed

Now we can state the announced main result. \\

\textbf{Theorem 1 (Riemannian Positive Mass Theorem) }{\itshape \, Let $(M^n,g)$ be asymptotically flat of order $\tau  > \frac{n-2}{2}$ with $scal(g)\ge 0$.\\

Then the total mass $E(M,g)$ is non-negative, $E \ge 0$, and $E = 0$ if and only if $(M,g)$ is isometric to the Euclidean space $\R^n$.}\\

We may also consider asymptotically flat manifolds with several such ends, but for the positive mass theorem there is a simple argument how to reduce the problem to the single end case. Namely, one deforms and closes all but one end to  large compact $scal>0$-bubbles, cf. for instance, [L1],Ch.5. In turn, this means, the theorem shows that the mass is non-negative for each of the ends.\\

Theorem 1, although physically motivated, obviously belongs to the field of scalar curvature geometry. The equality case $E = 0$ readily implies the following basic obstruction
result for positive scalar curvature metrics. It shows that there is no general mechanism to locally \emph{increase} the scalar curvature of arbitrary Riemannian manifolds even when we allow topological changes.\\

\textbf{Theorem 2} \textbf{(Non-Existence of $Scal > 0$-Islands) }{\itshape \,There exists \textbf{no} complete Riemannian manifold $(M^n, g)$ such that:
 \begin{itemize}
   \item $scal(g) >0$ on some non-empty open set $U \subset M^n$, with compact closure.
   \item $(M^n \setminus U, g)$ is isometric to $(\R^n \setminus  B_1(0), g_{Eucl})$.
 \end{itemize}}
\medskip
Less obviously, Theorem 2 also implies the apparently much stronger Theorem 1. This has been proved by the author in [L1],Ch.6.\\

 The idea to derive the non-negativity of $E$ from Theorem 2 is as follows: assume, in Theorem 1,  that $E < 0$, then we can bend the, therefore,  hyperbolically flavoured end to a flat Euclidean end. This compresses the volumina in a way that increases the scalar curvature and we get a $scal > 0$-island. \\

The rigidity statement, for $E = 0$, then follows from the non-negativity assertion for $E$ using  elementary arguments in conformal geometry. \\

\textbf{Remark 2 (Black Holes)} \, Both theorems seamlessly extend to cases where $M$ contains black holes. In the Riemannian situation, this means the non-asymptotically flat portion of $M$ contains (outer) area minimizing boundary components. We discuss this extension in \ref{stm} below.
 \qed

\subsubsection{Classical and New Techniques} \label{int2}
\bigskip

We discuss the limits of the classically known techniques to settle the positive mass conjectures. Then we recall some of the new skin structural methods from [L1]-[L3] we use to understand the general situation.\\

\textbf{Classical Approaches} \,  As mentioned above, the  positive mass conjectures had been established when either dimensions $\le 7$ or the underlying space is a spin  manifold.\\

The proofs, in dimensions $\le 7$, are owing  to Schoen and Yau [SY1]-[SY3] and [S], using minimal hypersurface techniques.  The reason why they only apply in dimensions $\le 7$ is that their argument only applies for \emph{smooth} area minimizing hypersurfaces. However, in dimensions beyond $7$, these hypersurfaces may carry delicate singularities and this causes some, by classical means, hardly resolvable issues, we discuss below.\\

The elegant spin theoretical approach to these problems, owing to Witten [W],[B] and [PT]  broadly apply in low dimensions, since, for instance, every orientable three manifold is spin. However, for increasing dimensions the \emph{non-spin} case gradually becomes the generic one, whereas the required spin condition turns into a restrictive symmetry-style constraint.\\

\textbf{Minimal Hypersurface Techniques} \, Our argument, to derive the general case of the Riemannian positive mass conjecture, merges  basic ideas of the minimal hypersurface approach with the recent skin structural techniques to understand the specific geometric analysis on singular minimal hypersurfaces, developed in [L3]-[L5]. To understand this strategy, we recall how minimal hypersurfaces are used in this context and discuss the role of their singularities.\\

The method of studying positive scalar curvature problems on three-manifolds by means of a \emph{dimensional descent} to minimal surfaces has been introduced by  Schoen and Yau [SY1]. Later, they also applied a variant of this technique to solve the three-dimensional positive mass conjecture [SY2], [SY3].\\

The key insight was, that stable minimal surfaces admit $scal >0$-metrics, after certain conformal deformations, when their ambient three-manifold has $scal >0$. \\

 This \emph{$scal >0$-heredity principle} still holds in higher dimensions and it suggested to also study higher dimensional scalar curvature problems and positive mass theorems, using an inductive dimensional descent along towers of nested area minimizing hypersurfaces with inherited curvature conditions until one reaches lower dimensional $scal>0$-geometries already properly understood.\\

\textbf{Issues and Results in Dimensions} $\ge 7$ \, But this program could only be accomplished in dimensions $< 7$ since the inductive step essentially employs the smoothness of area minimizing hypersurfaces, and this is granted only in these low dimensions. Here we reach the notorious complication, that obstructed extensions to higher dimensions: in dimensions $\ge 7$ such hypersurfaces $H^n$ generally contain delicate singular sets $\Sigma_H \subset H$. We describe a few of the resulting issues:
\begin{itemize}
  \item We need to understand certain classical elliptic operators, defined on $H \setminus \Sigma_H$, to control the deformation to a $scal>0$-metric. But this is rather difficult, since $H \setminus \Sigma_H$ and also these operators degenerate towards $\Sigma_H$ in an intricate way and there is no efficiently applicable structure theory for these singular sets.
  \item  In towers of inductive descent, which include singular area minimizers, we encounter some peculiar phenomena. The minimizers could (at least partially) disappear in the singular sets of their ambience. Also the minimizers can now become singular even in dimensions $\le 6$, since singularities of minimizers in singular spaces could have smaller codimension. This makes it difficult to draw non-trivial conclusions, e.g. since every surface admits $scal>0$-metrics with finitely many singularities.
\end{itemize}
But, and this is the basis for our approach, with more recent tools, in particular from geometric analysis on metric spaces, we can also employ $\Sigma_H$ to gain control. The theory of \emph{skin structures} on area minimizers, developed in [L3] and [L4] sheds light on the impact of $\Sigma_H$ on the geometry and analysis of, and also on, $H \setminus \Sigma_H$.
\begin{itemize}
  \item We get a \emph{hyperbolic unfolding} of $H \setminus \Sigma_H$: there is a canonical conformal deformation of  $H \setminus \Sigma_H$ to a complete, Gromov hyperbolic space of bounded geometry, its \emph{skin geometry}. And its \emph{Gromov boundary}  is homeomorphic to  $\Sigma_H$.
  \item   For broad classes of elliptic operators, the \emph{skin adapted operators}, this inherent hyperbolicity of  $H \setminus \Sigma_H$ can be exploited to derive controls  for their potential theory and asymptotic analysis, near $\Sigma_H$, on the original space $H \setminus \Sigma_H$. For instance, their \emph{Martin boundary} on $H \setminus \Sigma_H$ is  again  homeomorphic to $\Sigma_H$.
\end{itemize}
In [L5] we employed this hitherto unexpected degree of control over the elliptic analysis on $H \setminus \Sigma_H$ to establish the broader framework of a \emph{$scal >0$-heredity with surgery}.\\

  It allows us to implement strategies of inductive descent which include procedures to regularize all singularities to regular ends and closures while keeping $scal >0$.\\

In this paper, we take advantage of these new tools to derive Theorem 2. In the next section we prepare the setup and describe the main steps of the proof which follows in Ch.\ref{st}. We recall some basic concepts and results form [L3]-[L5] in Ch.\ref{sk} and Ch.\ref{sus}. We also recommend to have a glance at [L3],Ch.1.1 for a general introduction to this theme.\\

\subsubsection{The Scheme of Inductive Descent} \label{sk1}
\bigskip

We formulate a $scal >0$-heredity with surgeries suitable to derive Theorem 2. Theorem 1 then follows from our reduction argument in [L1],Ch.6, and some standard arguments, in Ch.\ref{bo}, for the case $E=0$.\\

\textbf{Overview} \, To prove  Theorem 2  we assume the existence of a hypothetical $scal > 0$-island. We transplant this island onto an actually very large torus and get, after some arbitrarily $C^3$-small perturbation,  a $scal >0$-manifold of the form $Q^n=T^n \cs N^n$, for some compact manifold $N^n$.\\

We get an area minimizing hypersurface $X^{n-1} \subset T^n \cs N^n$ homologous to $T^{n-1} \subset T^n \cs N^n$. We use the regularity theory for such hypersurfaces to infer that $X^{n-1}$ also contains a nearly flat and large torus component. That is, $X^{n-1}$ can again be written $X^{n-1}= T^{n-1} \cs N^{n-1}$, for some compact but now usually singular space $N^{n-1}$.\\

Here we reach the main step. By means of skin structural techniques we prove that for any neighborhood $V$ of the singular set $\Sigma_{N^{n-1}}$, we can conformally deform $X^{n-1}$ into some space $Z^{n-1}$  so that there is a smaller smoothly bounded neighborhood $U \subset V$ of $\Sigma_{Z^{n-1}}=\Sigma_{N^{n-1}}$ such that, relative  the new metric:
\begin{itemize}
\item $Q^{n-1}=Z^{n-1} \setminus U$ has $scal >0$ and $\p Q^{n-1}=\p U$ has positive mean curvature.
\end{itemize}
Our sign convention is that $S^{n-1} \subset \R^n$, viewed from inside, is positively mean curved.\\

 The condition for $\p Q^{n-1}$ means we build a \emph{shielding horizon} around $\Sigma_{Z^{n-1}}$. Area minimizing sequences of hypersurfaces in $Q^{n-1}$, homologous to $T^{n-2} \subset T^{n-1} \cs N^{n-1}$, will be repelled from $\p Q^{n-1}$ and we find an area  minimizing hypersurface $X^{n-2} \subset Z^{n-1} \setminus \overline{U}=Q^{n-1} \setminus  \p Q^{n-1}$. The simple but essential observation is that, from the viewpoint of $X^{n-2}$, its ambient space $Q^{n-1}$ is entirely smooth.\\

Finally, we show that the conformal deformation, that transformed $X^{n-1}$ into  $Z^{n-1}$, reproduces the second geometric main property of $T^n \cs N^n$:
 \begin{itemize}
\item $Q^{n-1}$ also contains a nearly flat and large torus component.
\end{itemize}

From this point onwards, we can proceed inductively and we eventually reach a smooth $scal >0$-surface $Q^2$ with a torus component. In particular, it has genus $\ge 1$.  But the Gau\ss-Bonnet theorem says that $Q^2$ cannot admit a $scal >0$-metric. This contradiction shows that there is no $scal > 0$-island.\\

Now we make the initial step of a compactification more explicit, formulate the details of the inductively repeated main step and properly infer Theorem 1. The proof of the inductive main step will occupy most of the remainder of this paper. Only in Ch.\ref{bo} we consider the rigidity case $E=0$, to complete the proof of Theorem 2. \\

 \textbf{Compactifications} \, We start  from a hypothetically existing  complete manifold $(M^n, g)$ with $scal(g) >0$ on some non-empty open set $U \subset M^n$, with compact closure $K$, so that $(M^n \setminus U, g)$ is isometric to $(\R^n \setminus  B_1(0), g_{Eucl})$.\\

$\bullet$ \, We consider the unit cube $C^n = [ -1 , 1] \times ...  \times [ -1 , 1] \subset \R^n$ and identify opposite sides of this cube to get the flat torus \[T^n= C^n/\sim=(S^1 \times ... S^1, g_{S^1} \times  ... \times  g_{S^1}),\]

where, in this case, $g_{S^1}$ is normalized to $length(S^1)=2$. We later consider the two distinguished points $1_T:=(1,...,1)/\sim$ and $0_T:=(0,...0)/\sim$. For any $x \in T^n$, we define the following
  $(n-1)$-dimensional torus that passes through $x$: \[T^{n-1}[x]:=\{x_1\} \times S^1 \times ... S^1 \subset T^n.\]

$\bullet$ \, We scale $(M,g)$ by some $\gamma \in (0,1/100)$, that is, we consider $(M^n, \gamma^2 \cdot g)$.\\

Now we insert $(K, \gamma^2 \cdot g) \subset (M^n, \gamma^2 \cdot g)$ into $T^n$: we delete the ball $\overline{B_\gamma(0)}$ from $T^n$   and replace it for $K$. This gives a smooth closed Riemannian manifold
$(Q^n_\gamma,g_{n,\gamma})$. Topologically, we have $Q^n=T^n \cs N^n$, for some closed manifold $N^n$ and we observe
\begin{enumerate}
  \item $Q^n_\gamma$ contains  a flat torus component $T^n\setminus \overline{B_\gamma(0)}=Q^n_\gamma \setminus K$.
  \item $scal(g_{n,\gamma}) \ge 0$ and $scal(g_{n,\gamma})> 0$ on the interior of $K$.
\end{enumerate}

A slight perturbation of this metric gives us $scal >0$-geometries  on $Q^n_\gamma$ we use as the starting case in our induction:
\begin{lemma} \emph{\textbf{(Base Case)}}\label{it}\,  For $n \ge 3$, any given $\ve, \gamma >0$, and $l \in \Z^{\ge 4}$, we can conformally deform $g_{n,\gamma}$ on $Q^n_\gamma$ into a smooth metric $g_n(\ve,\gamma,l)$ with
 \begin{equation}\label{bc}
 scal(g_n(\ve,\gamma,l))> 0 \mm{ on } Q^n_\gamma \,\mm{ and }\, |g_n(\ve,\gamma,l)-g_{n,\gamma}|_{C^l(Q^n_\gamma)}\le \ve.
 \end{equation}
\end{lemma}

\textbf{Proof} \quad The variational characterization of the first eigenvalue $\lambda_1$ of the conformal Laplacian for the original metric $g_{n,\gamma}$:  $L_{Q^n_\gamma}  f  = -\Delta  f +\frac{n-2}{4 (n-1)} \cdot scal_{Q^n_\gamma} \cdot f,$
{\small \[\lambda_1 = \inf \Big\{\int_{Q^n_\gamma} | \nabla f |^2 + \frac{n-2}{4 (n-1)}  \cdot  scal_{Q^n_\gamma}  \cdot  f^2 d A \, |\, f \mm{ smooth}, |f|_{L^2}=1 \Big\},\]}
shows  from (ii) above that $\lambda_1 > 0$.  Then, the transformation law for scalar curvature under conformal transformation by the first eigenfunction $f_1>0$, which by standard elliptic theory is a smooth function, shows
\begin{equation} \textstyle\label{4} scal(f_1^{4/n-2} \cdot g_{n,\gamma}) \cdot f_1^{\frac{n+2}{n-2}} = \gamma_n \cdot L_{Q^n_\gamma} f_1 =   \gamma_n \cdot \lambda^L_1 \cdot f_1 > 0. \mm{ for } \gamma_n = \frac{4 (n-1)}{n-2}.\end{equation}
Now we choose the smooth function $f[\ve]:= 1 + \ve  \cdot f_1$, for  $\ve >0$ and observe that
\begin{equation} \label{4a} scal(f[\ve]^{4/n-2} \cdot g_{n,\gamma}) \cdot f[\ve]^{\frac{n+2}{n-2}} =\end{equation}
 \[\textstyle \gamma_n \cdot L_{Q^n_\gamma} f[\ve]  =\gamma_n \cdot (\ve \cdot  \lambda_1 \cdot f_1 + \frac{n-2}{4 (n-1)} \cdot scal_{Q^n_\gamma}\cdot f[\ve]) > 0.\]
Also,  $f[\ve] \ra 1$ in $C^l$-norm, for $\ve \ra 0$. Therefore, for any $\ve >0$, we can find some $\ve^*(\ve)>0$, so that  $g_n(\ve,\gamma,l):=f[\ve^*]^{4/n-2} \cdot g_{n,\gamma}$ has the asserted properties. \qed

\textbf{Inductive Step} \,   Next we assert that, given such spaces $(Q^n_\gamma, g_n(\ve,\gamma,l))$, we can inductively find $scal >0$-manifolds with essentially the same properties, in all lower dimensions $m$, for $n > m \ge 2$.

\begin{proposition}\emph{\textbf{(Inductive Descent)}}\label{id} \, We assume that, for some $m > 2$ and any $\ve>0$, $\gamma \in (0,1/100)$, $l \in \Z^{\ge 4}$, there is a $C^l$-regular $m$-dimensional Riemannian manifold $(Q^m(\ve,\gamma,l),g_m(\ve,\gamma,l))$ with the following properties
\begin{enumerate}
\item $(Q^m,g_m)$ is a compact $scal >0$-manifold with boundary $\p Q^m$.
\item The boundary $\p Q^m$ may be empty and, for $m <7$, it actually always is empty.
\item  $Q^m$ contains a torus component $QT^m$, $\ve$-almost isometric to the flat $T^m\setminus B_{\gamma}(0)$:
\[\mm{There is a $C^l$-regular diffeomorphism } F_m(\ve,\gamma,l): QT^m \ra T^m\setminus  B_{\gamma}(0)\mm{ with }\] \[ |F_m^*(g_{flat})- g_m|_{C^l(QT^m)} \le \ve.\]
\item  When $\p Q^m \n$, then $\p Q^m  \subset Q^m \setminus QT^m$ is  $C^3$-smooth and compact with positive mean curvature.
\end{enumerate}
Then, there exist manifolds $(Q^{m-1}(\ve^*,\gamma^*,l^*),g_{m-1}(\ve^*\,\gamma^*,l^*))$ which
satisfy the $(m-1)$-dimensional version of these properties,  for any $\ve^*>0$,  $\gamma^* \in (0,1/100)$ and $l^* \in \Z^{\ge 4}$.
\end{proposition}

One can think of this scheme as reducing the proof of the non-existence of  $scal > 0$-islands in dimension $m$ to similar non-existence results in lower dimensions.\\

The surface case $m-1=2$ is special, in the sense that the used deformations follow a slightly different pattern. To avoid a separate treatment of this case, we only show that $Q^2$ has genus $\ge 1$ and a positive integral scalar curvature. This suffices to derive a contradiction from the Gau\ss -Bonnet theorem, whereas all other conditions are only used in the induction hypotheses of the, now no more existing,  next step.

\begin{remark}\, We briefly explain the origin and use of some of the details of \ref{id}:\\

The manifolds $Q^{m-1}$ will be deformed and regularized minimal hypersurfaces in the $Q^m$. The non-trivial boundary components are offshoots of regularized singularities. We get these boundaries from a $scal >0$-heredity with surgery.\\

$\bullet$\,  In the low dimensional cases, without singularities, we do not need any surgeries and we can track the inherited topology from $Q^m$ to $Q^{m-1}$, by purely topological means from iterative use of duality arguments between integral homology and de Rham cohomology. \\

For $n< 7$, we start with a degree $1$ map $Q^n \ra T^n$, to find linear independent closed $1$-forms $\theta_i$, $i=1,..,n$ on $Q^n$, so that $T^{n-1}[1_T] \subset Q^n$ is dual to $\theta_1 \wedge.. \theta_{n-1}$ and $\int_{T^{n-1}[1_T]}\theta_1 \wedge.. \theta_{n-1}=1$. We can iteratively find a tower of smooth hypersurfaces, each conformal to a $scal >0$-manifold and we  finally reach a surface $F$ with $scal >0$ so that $\int_F \theta_1 \wedge \theta_2=1$, for two closed $1$-forms $\theta_i$, $i=1,2$, on $F$. But then the surface must have genus $\ge 1$, contradicting the Gau\ss -Bonnet theorem. \\

In the general case, $n \ge 7$, we have $\Sigma \n$ and iteratively apply $scal >0$-heredities with surgeries. This makes it hard to employ this topological tracking. Instead we use the robust geometric constraint of an almost flat torus component in $Q^m$, which survives under dimensional descent, even for $\Sigma \n$.\\

$\bullet$\,  The rather tecchnical $C^l$- and $C^{l^*}$-regularity conditions are introduced to compensate the decreasing regularity we observe during the process of inductive descent. Recall that minimal hypersurfaces within a $C^k$-regular manifold will only be $C^{k-1}$-regular, outside their singular set [P],Th.7.2 and [F],5.3.20. \\

This causes an issue not resolvable from choosing $l \ge n+2$, since the plain $scal >0$-heredity with surgeries, of [L5], reduces the regularity even of a smooth hypersurface  to  $C^{2,\alpha}$-regularity.  Therefore, we inductively apply a smoothing step, in Ch.\ref{sus}, with reference to the given regularity, to keep the regularity high enough to reach a $C^2$-surface.\\

We choose $Q^n(\ve,\gamma,l)=Q^n_\gamma$, for any $\ve, l$, whereas the topology of the lower dimensional $Q^m(\ve,\gamma,l)$, for $m <n$, may  depend on both $\ve$ and $l$.\qed
\end{remark}

 From Lemma \ref{it} and Prop.\ref{id}, we can easily derive our Theorem 2 and from this we infer, using also Ch.\ref{bo}, Theorem 1.

\begin{corollary} \emph{\textbf{(Theorem 2)}}  \,\label{teo} There are no $scal > 0$-islands.
\end{corollary}

\textbf{Proof} \quad The family of manifolds $(Q^n(\ve,\gamma,l), g_n(\ve,\gamma,l)):=(Q^n_\gamma, g_n(\ve,\gamma,l))$ satisfies the assumptions of \ref{id} in dimension $n$. Thus we may iteratively apply the inductive descent of \ref{id}, keeping an inheritable torus component until we reach a surface of genus $\ge 1$ with positive (integral) scalar curvature. But the Gau\ss -Bonnet theorem excludes the existence of such a surface. This contradiction shows that there are no $scal > 0$-islands.\qed

This makes the \textbf{proof} of Prop.\ref{id} our man objective. We describe the main steps. For the present we omit the  technical $C^l$-regularity considerations.
\begin{itemize}
\item  Ch.\ref{ami}: We use the regularity theory of area minimizing hypersurfaces to get a, potentially singular, area minimizer $X^{m-1} \subset (Q^m(\ve,\gamma),g_m(\ve,\gamma))$, so that $X^{m-1}$ contains a torus component $XT^{m-1}$, which is $\ve^*$-almost isometric to the flat $T^{m-1}\setminus \overline{B_{\gamma^*}(0)}$, with $\gamma^* \ra 0$, for $\gamma \ra 0$ and $\ve^* \ra 0$, for $\ve \ra 0$.
 \item  Ch.\ref{sk}: We recall some skin structural concepts, from [L3] and [L4], we need to define and to control conformal deformations in the next two sections.
\item  Ch.\ref{sus}: Now we conformally deform $X^{m-1}$ to a $scal >0$-manifold $Y^{m-1}$. The conformal factor is an eigenfunction of a skin variant of the conformal Laplacian. We employ the particular analysis of these eigenfunctions in the final deformation step to close the inductive loop, in the next section.
 \item Ch.\ref{mos}: For $\Sigma_X \n$ we further conformally deform $Y^{m-1}$, near $\Sigma_X$, to a space $Z^{m-1}$, so that on the complement of a small neighborhood $U_Z$ of $\Sigma_Z$,  $Z^{m-1}$ is a $scal >0$-manifold $(Q^{m-1}(\ve^+,\gamma^+),g_{m-1}(\ve^+,\gamma^+))$ with a positively mean curved boundary  $\p U_Z$ and $Q^{m-1}$ now contains a component, $\ve^+$-almost isometric to the flat $T^{m-1}\setminus \overline{B_{\gamma^+}(0)}$, with  $\gamma^+ \ra 0$, for $\gamma^* \ra 0$,  $\ve^+ \ra 0$, for $\ve^* \ra 0$.  \qed
\end{itemize}

\setcounter{section}{2}
\renewcommand{\thesubsection}{\thesection}
\subsection{$Scal >0$-heredity with Surgeries and Constraints} \label{st}
\bigskip

\subsubsection{Geometry of Area Minimizers in $Q^m$} \label{ami}
\bigskip

Given the spaces $Q^m$ specified in \ref{id}, the first step to get the hypersurfaces  $Q^{m-1}$ is to choose some area minimizing hypersurfaces in  $Q^m$. This uses some  geometric measure theory. Details and explanations concerning these basics can be found in [F], [GMS], [KP], [P] and [Si1]. An overview of these results is in [L3], Appendix. For the reader's convenience and to fix our notations we recall some of the concepts and results.\\

$\bullet$ \,  \textbf{Currents}\,  The extension of the space of submanifolds to that of \emph{$m$-currents} $\mathcal{D}_m(U)$, with $U \subset \R^n$ open, is the dual of $\mathcal{D}^m(U)$, the space of smooth $m$-forms compactly supported in $U$. When we integrate such  forms over a submanifold $N$, we can interpret $N$ as a current, denoted by $\llb N \rrb$. \\

The weighted area of   $T  \in \mathcal{D}_{m}(U)$, its \emph{mass} ${\bf{M}}_U (T)$, is defined by
\[{\bf{M}}_U (T) = \sup \{T(\omega)\,|\, \omega \in \mathcal{D}_m(U), |\omega| \le 1, \mm{ supp }\omega \subset U\}.\]

A subset of $\R^n$ is called\emph{ k-rectifiable} if it is the image of a bounded subset of $\R^k$ under a Lipschitz map. A current associated to a rectifiable set is a \emph{rectifiable current}, where we allow integer multiplicities. If the boundary of such a current, defined using Stokes' formula, is also rectifiable it is called \emph{integral current}. The restriction of a rectifiable $T \in \mathcal{D}_m(U)$ to any Borel set $V \subset U$ is denoted $T \llcorner V$.\\

$\bullet$ \,\textbf{Oriented Boundaries} \, Any mass minimizing integral current $X$ locally, e.g.  in some ball $B$, is a union of minimal \emph{oriented boundaries} $X_i$ of open sets $U_i \subset B$, that is, $\p U_i \cap B= X_i \cap B$, [KP],7.5, [Si1],Ch.26. The regularity theory asserts that each $X_i$ is a smooth hypersurface except for singularities of at most $m-7$ dimensions, [KP],9, [Si1],37.\\

In general the union is countable, but in our case of a compact ambient space, it is a \emph{finite union}. Also, the strict maximum principle shows that any two $X_i$ are locally disjoint or they are identical, in cases of higher multiplicities. For simplicity, the global path-component of each $X_i$ in $X$, is still denoted by $X_i$. Globally, $X_i$ usually no longer represents an oriented boundary, since boundaries are null-homologous. Thus, we henceforth call it a \emph{minimal sheet} of $X$.\\

$\bullet$ \,\textbf{Flat Metric} \, For any open subsets of $\R^n$: $W \subset \overline{W} \subset U
\subset \mathbb{R}^n$, we define the \emph{flat (pseudo)metric} on $\mathcal{D}_{n-1}(U)$:
{\small
\[\db_W(C_1,C_2) := \mbox{inf}\{{\bf{M}}_W (S) +
{\bf{M}}_W (R)\, |\, C_1 - C_2 = S + \p R, S \in \mathcal{D}_{n-1}(U), R \in \mathcal{D}_n(U)) \}.\]}
for $C_1, C_2 \subset \mathcal{D}_{n-1}(U)$. In simple terms, it measures the volume between  $C_1$ and $C_2$.   The family of these $\db_W$ generate the flat metric topology. Using charts these ideas readily extend to manifolds. \\

Now we reach the basic existence result we need in our further discussion:

\begin{proposition}\emph{\textbf{(Homological Minimizers)}}\label{hm} \, For any $\ve>0$, $\gamma \in (0,1/100)$ and $l \in \Z^{\ge 4}$, there is a mass minimizing integral current $X^{m-1}(\ve,\gamma,l-1)$ so that:
\begin{itemize}
\item  $X^{m-1} \subset (Q^m(\ve,\gamma,l),g_m(\ve,\gamma,l))$ and $X^{m-1} \cap \p Q^m \v$.
  \item  $X^{m-1}$ is homologous to the torus $F_m^{-1}(T^{m-1}[1_T]) \subset QT^m$.
  \item The support of  $X^{m-1}$ is a $C^{l-1}$-regular hypersurface, except for some singular set $\Sigma_X$ of codim $\ge 8$ in $Q^m$.
\end{itemize}
\end{proposition}

\textbf{Proof} \quad This largely is a classical result, we composed from  [GMS], 5.4.1,Cor.1 and [P],Th.7.2. The only, however, important detail we need to explain is that, this holds even when $\p Q^m \n$ and $X^{m-1} \cap \p Q^m \v$.\\

Let $H_k$, $k \in \Z^{\ge 1}$, be a mass minimizing sequence of hypersurfaces in $(Q^m,g_m)$ homologous to $F_m^{-1}(T^{m-1}[1_T]) \subset QT^m$. Since $\p Q^m  \subset Q^m \setminus QT^m$ has positive mean curvature it repels a mass minimizing sequence of hypersurfaces, that is, when
a sequence $\p Q^m$ tends to approach it, there is a competing sequence $H^*_k$ with smaller mass supported in positively bounded distance from $\p Q^m$.\\

 The basic compactness result for such currents, [KP],Th.7.5.2,  gives us a converging subsequence of $H^*_k$ with limit minimizer $H_\infty$  homologous to $F_m^{-1}(T^{m-1}[1_T])$ in $(Q^m,g_m)$ and we have $H_\infty \cap \p Q^m \v$.  \qed

\begin{remark}\textbf{(Minimal Boundary Components)}\label{stm} \, 1. We notice that the repelling effect of $\p Q^m$ also shows that there is an area minimizing hypersurface $L_m \subset Q_m$ homologous to $\p Q^m$, that surrounds
$\p Q^m$. That is, $L_m$ separates $Q^m$ in two components $Q_i^m$, $i=1,2$, with $\p Q^m \subset Q_1^m$ and $QT^m \subset Q_2^m$. Moreover, we can choose $L_m$ to be the hypersurface closest to $\p Q^m$.\\

Then the limit minimizer $H_\infty$ does not intersect $Q_1^m$ and, due to the strict maximum principle [Si2], we even have $H_\infty \cap L_m \v$, regardless the regularity of $L_m$.\\

2. In the introduction we noticed, that both Theorems still apply to spaces whose non-asymptotically flat portion of $M$ contains (outer) area minimizing boundary components. The latter argument can also be used to include this black hole case. The one modification in \ref{it} needed is to use a conformal deformation $g_{n,\gamma}$ on $Q^n_\gamma$ again to  a smooth metric $g^*_n(\ve,\gamma,l)$ with
\[scal(g^*_n(\ve,\gamma,l))> 0 \mm{ on } Q^n_\gamma \,\mm{ and }\, |g^*_n(\ve,\gamma,l)-g_{n,\gamma}|_{C^l(Q^n_\gamma)}\le \ve,\]
 with the additional property that $g^*_n(\ve,\gamma,l) \equiv g_{n,\gamma}$ near $\p Q^m$. This can be accomplished by elementary means, similar to  [L2],Prop.1.1. \qed
\end{remark}

Next we further exploit the finer assumptions on the properties of $Q^m$ to sharpen the geometric information we may get from the regularity theory.\\

For given $\gamma >0$ and $\gamma^*\ge \gamma$, we henceforth write
\begin{equation}\label{ga}
QT^m_{\gamma^*}:= F_m(\ve,\gamma,l)^{-1}(T^m\setminus  B_{\gamma^*}(0)) \subset QT^m (= QT^m_{\gamma})
\end{equation}

\begin{proposition}\emph{\textbf{(Regularizing Symmetry)}}\label{rs} \, For any $\ve^* >0$, $\gamma^* \in (0,1/100)$ and $l^* \in \Z^{\ge 4}$, there are some $\ve \in (0,\ve^*)$, $\gamma \in (0,\gamma^*)$ and   $l = l^*+1$, so that there is a current $X^{m-1} \subset (Q^m(\ve,\gamma,l), g_m(\ve,\gamma,l))$, as in \ref{hm}, with the following geometric properties:
\begin{enumerate}
  \item $X \cap  QT^m_{\gamma^*}$ is a $C^{l^*}$-regular hypersurface.
  \item  $X \cap  QT^m_{\gamma^*}$  is almost isometric to $T^{m-1}[0_T] \setminus B_{\gamma^*}(0_T)$:
  \begin{equation}\label{f0}
 X \cap   QT^m_{\gamma^*}= exp_\nu(\Gamma)  \cap QT^m_{\gamma^*},
 \end{equation}
for  the image $exp_\nu(\Gamma)$ of some $C^{l^*}$-regular section $\Gamma$ of the normal bundle $\nu$ over $F^{-1}_m(T^{m-1}[0_T] \setminus   B_{\gamma^*}(0_T))$,  under the exponential map $exp_\nu$ of $\nu$, with
 \begin{equation}\label{f}
 |\Gamma|_{C^{l^*}\left(F^{-1}_m(T^{m-1}[0_T] \setminus   B_{\gamma^*}(0_T))\right)} <\ve^*.
 \end{equation}
\end{enumerate}
\end{proposition}

\textbf{Proof} \quad Obviously, (ii) implies (i), since it gives a presentation of $X \cap   F_m^{-1}(T^m \setminus B_{\gamma^*}(0_T))$ as a $C^{l^*}$-regular section. \\

 We now turn to the proof of (ii). The argument essentially is a stable characterization of the subtori of a flat torus as the unique area minimizers in their homology class.\\

 We first consider the case of an entirely flat torus $T^m$, we think of it as (the completion of the) limit of
$T^m  \setminus B_\gamma(0)$, for $\gamma \ra 0$. Then we treat  the still flat case  $\ve =0$, $\gamma >0$ and finally the general case $\ve >0$, $\gamma >0$ as two consecutive perturbations.\\

\textbf{1. Flat Tori} ($\ve =0, \gamma=0$)\, In the flat $T^m$, any area minimizer $X(0,0) \subset T^m$ homologous to  $T^{m-1}[1_T]$ is also a flat subtorus, $X(0,0)=\llb T^{m-1}[p]\rrb$, for some $p \in T^m$. Since we could not localize this folklore result in the literature we include a short proof. \\

We use the  constancy theorem, [KP],7.3, [Si1],Th.26.27.  It asserts that for  any $T \in \mathcal{D}_m(U)$, $U \subset \R^m$ open, connected and $\p T \v$ relative $U$, there is some $c \in \R$ such that $T=c \cdot \llb U \rrb$, and we have $c \in \Z$, when $ T$ is an integral current.\\

For the projection $\pi: T^m \ra T^{m-1}[1_T]$, this shows $\pi_*(X_i)= c_i \cdot\llb T^{m-1}[1_T]\rrb$, for some $c_i \in \Z$. The minimality of each minimal sheet $X_i \subset X$ and of $\llb T^{m-1}[1_T]\rrb$ then implies that there is only one non-trivial sheet, let us say $X_1$, and  $c_1=  1$. In particular, we have $S^1 \times \{a_2\} \times ... \{a_n\} \cap X_1 \n$, for any
$(a_2,..a_n) \in T^{m-1}[0_T]$.\\

Now assume that $X_1$ is not contained in some $T^{m-1}[p]$, then, since its singularities have at most $m-8$ dimensions, it intersects $T^{m-1}[p]$, for a suitable $p \in T^m$, transversely along some $(m-1)$-dimensional submanifold, in some open subset of $T^m$ disjoint to $\Sigma_{X_1}$. From this we get another minimal sheet homologous to  $T^{m-1}[0_T]$ with an $(m-1)$-dimensional singular set.  This contradicts the fact that the singularities (of the new minimizer) can be at most $m-8$-dimensional and we infer $X(0,0)=\llb X_1\rrb=\llb T^{m-1}[p]\rrb$, for some $p \in T^m$.\\

\textbf{2. Flat Components} ($\ve =0, \gamma >0$)\, We show that for  $\gamma >0$, and any mass minimizer  $X(0,\gamma,l^*) \subset Q^m(0,\gamma,l^*+1)$ homologous to  $T^{m-1}[1_T]$, we have a yet to specify type of (sub)convergence of pairs
\begin{equation}\label{kin}
 (X(0,\gamma,l^*), Q^m(0,\gamma,l^*+1)) \, \ra \, (T^{m-1}[p],  T^m), \mm{ when } \gamma \ra 0
\end{equation}
$X(0,\gamma,l^*) \subset Q^m(0,\gamma,l^*+1)$,  $T^{m-1}[p]\subset  T^m$, for some $p \in T^m$.\\

We may assume that for each minimal sheet $X_i \subset X$: $X_i \cap T^m \setminus QT^m \n$. Otherwise, we are in Case 1 and there is a non-trivial minimal
sheet $X_i \subset QT^m$. Then the constancy theorem and the minimality of $X_i$ and $\llb T^{m-1}[1_T]\rrb$ show again that $X=X_i$ is a flat subtorus.\\

We plan to apply an compactness argument to the $ X(0,\gamma)$, for $\gamma \ra 0$ in (\ref{kin}), but there is no suitable control over the non-torus components $Q^m(0,\gamma) \setminus QT^m(0,\gamma)$, under variations of $\gamma$. Therefore, we consider the restrictions  $X(0, \gamma) \llcorner QT^m$. Indeed, we henceforth only use that
 \begin{itemize}
   \item $X(0, \gamma) \llcorner QT^m$ is area minimizing for perturbation compactly supported in $QT^m$
   \item $(\p X(0, \gamma)) \llcorner QT^m \v$ and $X(0, \gamma) \llcorner QT^m$  is (relative) homologous to $\llb T^{m-1}[1_T]\rrb$.
 \end{itemize}

Now we get mass estimates from comparisons with combinations of elementary currents.
\begin{equation}\label{m1}
{\bf{M}}( X(0,\gamma)\llcorner T^m \setminus B_{2 \cdot \gamma}(0)) \le {\bf{M}}(T^{m-1}[1_T]) + {\bf{M}}(S^m).
\end{equation}
Also, because $X(0,\gamma)$ is still area minimizing relative $T^m_\gamma\setminus B_\gamma(0_T)$,  the coarea formula, [KP],5.2, and specifically [GMS],2.1.5, Th.3, and the slicing theorem, [KP],7.6, [Si1],13,  say that for each $\gamma$ there is some radius $\rho_\gamma \in (\gamma,2 \cdot \gamma)$ so that $\p ( X(0,\gamma)\llcorner  B_{\rho_\gamma}(0_T))$ is a rectifiable current and
 \begin{equation}\label{m2}
 {\bf{M}}\left(\p ( X(0,\gamma)\llcorner  B_{\rho_\gamma}(0_T))\right) \le c_m, \mm{ for some } c_m>0, \mm{ depending only on  } m.
 \end{equation}
Since the $ X(0,\gamma) \llcorner T^m\setminus B_{\rho_\gamma}(0_T)$ are supported in $T^m$ and satisfy the uniform bounds (\ref{m1}) and (\ref{m2})
 we may apply  the compactness theorem for integral currents, [KP],Th.7.5.2, [Si1],Th.32.2, and find, for any sequence $\gamma_i \ra 0$, for $i \ra \infty$, a $\db$-converging subsequence of the $X(0,\gamma_i) \llcorner T^m\setminus B_{\rho_{\gamma_i}}(0_T)$ on $T^m$. For any  limit $P \subset T^m$, the semi-continuity of the area (mass) functional [KP],(7.5),p.187, shows that $P$ is mass minimizing relative any perturbation compactly supported in $T^m\setminus \{0\}$ and homologous to $\llb T^{m-1}[1_T]\rrb$.\\

Also, we observe for the boundary of the tube $S^1 \times B_r(0)\subset S^1 \times T^{m-1}[0_T]$:  \[{\bf{M}}(P) \le {\bf{M}}(\llb T^{m-1}[0_T]\rrb)  + {\bf{M}}(\llb  \p (S^1 \times B_r(0))\rrb) \mm{ and }{\bf{M}}(\llb \p (S^1 \times B_r(0)) \rrb) \ra 0, \mm{ for }r \ra 0.\] From this we infer that $P$ is mass minimizing in $T^m$.\\

Hence, from Case 1 we conclude that $P=\llb T^{m-1}[0_T]\rrb$. From this we also see that actually the entire sequence $X(0,\gamma_i) \llcorner T^m\setminus B_\beta(0_T)$  $\db$-converges to $\llb T^{m-1}[0_T]\rrb \llcorner T^m\setminus B_\beta(0_T)$, for any $\beta>0$ given.\\

\textbf{Graphical Representation} \, As a consequence of the monotonicity formula, [Si1],17,  there is a small $\ve_{m,l}>0$, so that for $\R^m$ equipped with a Riemannian metric $g$ with $|g-g_{Eucl}|_{C^{l^*}(B_2(0))}\le \ve_m$ and any minimal boundary $\p U \subset (\R^m,g)$, with $0 \in \p U$,  we have
\begin{equation}\label{aes}
Vol_{m-1}(\p U \cap B_1(0)) \ge \iota_m,\, \mm{ for some } \iota_m >0, \mm{ depending only on } m.
\end{equation}
Here we covered also the not strictly flat case as a preparation for the next step 3. \\

For the $\db$-converging sequence of area minimizers $X(0,\gamma_i)$ this shows that $ X(0,\gamma_i)$ compactly Hausdorff converges in $T^m\setminus \{0\}$ to $\llb T^{m-1}[0_T]\rrb$. Now Allard regularity theory for area minimizing hypersurfaces, [Si1],23 and 24,  shows that for any $\beta>0$, there is an $i_\beta$, so that
\begin{description}
  \item[G1] For $i \ge i_\beta$,  $ X(0,\gamma_i) \llcorner T^m \setminus S^1 \times B_\beta(0_T)$ can be represented by a  smooth graph $\Gamma_i$ over $T^{m-1}[0_T]$ for some smooth function $G_i$, where we locally think of $S^1$ as $\R$.
  \item[G2] For $i \ra \infty$ we have: $|G_i|_{C^k(T^{m-1}[0_T] \setminus B_{\beta}(0_T))} \ra 0$, for any $k \in \Z$.\\
\end{description}

\textbf{3. The General Case} ($\ve >0, \gamma>0$)\, Now we reach the case where $\ve >0$.
We note from the argument of Case 2, that for any area minimizer $X(\ve,\gamma,l^*)$ in $(Q^m(\ve,\gamma,l^*+1),g_m(\ve,\gamma,l^*+1))$ homologous to $\llb F_m^*(T^{m-1}[1_T]) \rrb$  we have a subconvergence
\[X(\ve,\gamma) \llcorner T^m\setminus B_{2 \cdot \gamma}(0_T) \ra X(0,\gamma) \llcorner T^m\setminus B_{2 \cdot \gamma}(0_T), \mm{ for } \ve \ra 0,\]
in $\db$-topology and thus, from (\ref{aes}), in Hausdorff-topology.\\

Then, another application of regularity theory for area minimizers shows that for any $\beta >0$ there are small $\ve, \gamma>0$ so that, from G1 and G2 above, $X(\ve,\gamma) \llcorner T^m\setminus B_{\beta}(0_T)$ is a $C^{l^*}$-smooth hypersurface arbitrary $C^{l^*}$-close to  $X(0,\gamma) \llcorner T^m\setminus B_{\beta}(0_T)$.\\

In turn, this means:
\begin{description}
  \item[G1*] For sufficiently small $\ve >0, \gamma>0$: $X(\ve,\gamma) \llcorner F^{-1}(T^m \setminus S^1 \times B_\beta(0_T))$ can be represented by the $F_m$-preimage of a $C^{l^*}$-regular graph $\Gamma_{\ve,\gamma}$ over $T^{m-1}[0_T] \setminus B_\beta(0_T)$ of a $C^{l^*}$-function $G_{\ve,\gamma}$.
  \item[G2*] For any $\eta >0$, we find small $\ve >0, \gamma>0$ so that $|G_{\ve,\gamma}|_{C^{l^*}(T^{m-1}[0_T] \setminus B_{\beta}(0_T))}\le \eta$.
\end{description}
This readily implies the assertions made in (ii).\qed

\subsubsection{Skin Structural Concepts} \label{sk}
\bigskip

In this section we describe some basic geometric analysis on singular area minimizing hypersurfaces by means of skin structures. We apply these details to define and control our subsequent conformal deformations of such hypersurfaces in the next section. We recollect some essentials from [L3] and [L4]. \\

\textbf{Skin Structures} \, We consider the class of area minimizing hypersurfaces $H^k \subset W^{k+1}$, where $W^{k+1}$ is a  compact manifold or the flat $\R^{k+1}$, for $k\ge 2$. For the sake of simplicity, we assume that $W^{k+1}$ and hence $H^k$ are orientable.\\

To any such hypersurface we assign a non-negative locally Lipschitz function $\bp_H$ on $H \setminus \Sigma$. The assignment
$H \mapsto \bp_H$ is a \emph{skin transform} provided it satisfies the following axioms:
\begin{itemize}
\item $H \mapsto \bp_H$  is naturally defined: it
commutes with the $\db$-convergence of sequences of area minimizers $H_i$, $i \ge 1$ to a
        limit space $H_\infty$, so that   \[\bp_{H_i} \overset{C^\alpha}  \longrightarrow {\bp_{H_\infty}}, \mm{ for any } \alpha \in (0,1)\]
  \item $\bp_H \ge |A_H|$, where $|A_H|$ is the norm of the second fundamental form of $H$.
        \item For any $f \in C^\infty(H \setminus \Sigma,\R)$ compactly supported in $H \setminus \Sigma$ and some $\tau = \tau(\bp,H) \in (0,1)$ we have the Hardy type inequality
\begin{equation}\label{hi}
\int_H|\nabla f|^2  + |A_H|^2 \cdot f^2 dA \ge \tau \cdot \int_H \bp_H^2\cdot f^2 dA.
\end{equation}
    \item $\bp_H \equiv 0$, if $H \subset N$ is totally geodesic, otherwise, $\bp_H$ is  strictly positive.
    \item When $H$ is not totally geodesic, we define $\delta_{\bp}:=1/\bp$, the \emph{${\bp}$-distance}. It is $L_{\bp}$-Lipschitz regular, for some constant
        $L_{\bp}=L(\bp,n)>0$:
        \[|\delta_{\bp}(p)- \delta_{\bp}(q)|   \le L_{\bp} \cdot d(p,q), \mm{ for } p,q \in  H \setminus \Sigma. \]
\end{itemize}
For a discussion and constructions of skin transforms, also explaining the terminology, we refer to  [L3],Def.1 and Ch.2.\\

\textbf{Skin Uniformity and Adaptedness} \, Skin transforms are used to express but also to prove otherwise hardly accessible geometric and analytic results for (and on) $H \setminus \Sigma$,  where we deliberately also allow the regular case $\Sigma \v$:\\

We get that $H \setminus \Sigma$ is a \emph{skin uniform space}, [L3], Th.2. That is, for any $p,q \in H \setminus \Sigma$, there is a rectifiable path $\gamma_{p,q}: [a,b] \ra H \setminus \Sigma$, for some $a <b$,
so that for any given skin transform $\bp$, there is some $s_H \ge 1$ with
\begin{enumerate}
\item \quad  $l(\gamma_{p,q})  \le s_H \cdot d_{g_H}(p,q)$, for any two points $p,q \in H \setminus \Sigma$
    \item \quad $l_{min}(\gamma_{p,q}(z)) \le s_H \cdot \delta_{\bp}(z)$, for any $z \in \gamma_{p,q}.$
\end{enumerate}
$l_{min}(\gamma_{p,q}(z))$ denotes the minimum of the lengths of the subarcs of $\gamma_{p,q}$ from $p$ to $z$ and from $q$ to $z$.\\

These conditions describe a form of boundary regularity of $\Sigma$ from $H \setminus \Sigma$, when we think $\Sigma$ as a boundary of $H \setminus \Sigma$. More generally, in the case where $\Sigma \v$, skin uniformity describes a quantitative approachability of subsets of $H$.  This uniformity is an essential ingredient in the study of \emph{skin adapted operators} $L$ over $H \setminus \Sigma$, [L4],Def.1, we may  characterize as follows:
\begin{enumerate}
\item $L$ is a linear second order elliptic operator and its coefficients do not diverge faster than $\bp^2$, while we approach $\Sigma$.
\item $L$ satisfies a weak coercivity condition:  there is a supersolution $s>0$ of $L \, u=0$ and some $\delta >0$, so that: $L \, s \ge \delta \cdot \bp^2 \cdot s.$
\end{enumerate}

For this quite large class of operators, the skin uniformity of $H \setminus \Sigma$ grants a transparent potential theory and asymptotic analysis near
$\Sigma$ described in [L4].\\

 This analysis is substantially employed to establish the results we discuss, and mostly cite from [L4], in the next two sections.\\

\textbf{Conformal Laplacians} \,  In this paper we are interested in operators related to scalar curvature geometry. The relevant operators are all variants of the conformal Laplacian $L_H:=-\Delta  +\frac{k-2}{4 (k-1)} \cdot scal_H$, for $k \ge 3$. (The simple but slightly differing surface case, $k=2$, is treated separately.)\\

This is not surprising since $L_H$ belongs to the $scal$- transformation law under conformal deformations $g_H \mapsto u^{4/k-2} \cdot g_H$:
\begin{equation}\label{stl}
\textstyle scal(u^{4/k-2} \cdot g_H) \cdot u^{\frac{k+2}{k-2}} =  \gamma_k \cdot L_H\,  u, \mm{ for } k \ge 3,   \mm{ and } \gamma_k = \frac{4 (k-1)}{k-2}
\end{equation}
valid for any $C^2$-regular function $u>0$. The simple but vital observation is that $L_H$ matches the skin structural techniques in the case of non-negative scalar curvature ambient spaces.

\begin{proposition}\emph{\textbf{(Conformal Laplacian $L_H$)}}\label{lh0} \, For $k \ge 3$, the conformal Laplacian $L_H$  is a skin adapted operator on $H \setminus \Sigma$, if and, in general, only if $scal_W \ge 0.$
\end{proposition}
In this case (of a symmetric operator) the weak coercivity requirement (ii) of the skin adaptedness conditions is equivalent to
\begin{equation}\label{cf}
\int_H  f  \cdot  L_H f  \,  dA  \ge \tau \cdot \int_H \bp^2 \cdot f^2   \, d A, \mm{ for some }\tau >0,
\end{equation}
for any smooth function $f$ on $H$ with $supp \: f \subset H \setminus \Sigma$.\\

\textbf{Proof of \ref{lh0}} \quad For later reference, we summarize the arguments from [L4], Th.10  and Prop.5.10. Since $\bp \ge |A|$
the Gau\ss-Codazzi equation, $|A|^2 + 2 \cdot Ric_W(\nu,\nu)  =  scal_W - scal_H$, readily shows that $L_H$ satisfies skin adaptedness condition (i).\\

For (ii), we recall that any infinitesimal variation normal to $H$, can be written as a normal field $f \cdot \nu$, where $\nu$ is the outward
normal vector field over $H \setminus \Sigma$, $f$ is a smooth function on $H$ with $supp \: f \subset H \setminus \Sigma$. The area minimizing property of $H$ shows that the the second variation of the area functional is non-negative: {\small
\begin{equation}\label{stabi}
Area''(f) = \int_{H}|\nabla_H f|^2 - \left( |A|^2 + Ric(\nu,\nu) \right) \cdot f^2 \: dVol \ge 0
\end{equation}
}
The Gau\ss-Codazzi equation gives an equivalent formulation of $Area'' (f) \ge 0$: {\small
\begin{equation} \label{kwsy}  \int_H  f  \cdot  L_H f  \,  dA  = \int_H | \nabla f |^2 + \frac{k-2}{4 (k-1)} \cdot scal_H  \cdot  f^2 \, d A \ge \end{equation}
\[\int_H \frac{k}{2 (k-1)} \cdot  |  \nabla f |^2 + \frac{k- 2}{4 (k-1)}\cdot \left( | A |^2 + scal_W \right) \cdot  f^2\,  d A\]}
Now we also assume that $scal_W \ge 0$, then we get from the Hardy inequality (\ref{hi}): {\small
\begin{equation} \label{aaa}  \int_H  f  \cdot  L_H f  \,  dA  \ge \int_H \frac{k}{2 (k-1)} \cdot  |  \nabla f |^2 + \frac{k- 2}{4 (k-1)}\cdot  | A |^2  \cdot  f^2\,  d A \end{equation}
 \[ \ge \tau(\bp,H)  \cdot \frac{k- 2}{4 (k-1)}\cdot \int_H \bp^2 \cdot f^2   \, d A.\]} \qed

 There is a largest $\lambda \ge \tau(\bp,H)  \cdot \frac{k- 2}{4 (k-1)}\tau$, so that $\int_H  f  \cdot  L_H f  \,  dA  \ge \lambda\cdot \int_H \bp^2 \cdot f^2   \, d A,$
for any smooth function $f$ on $H$ with $supp \: f \subset H \setminus \Sigma$.\\

Its value $\lambda^{\bp}_{H} >0$ is called the \emph{generalized principal eigenvalue} of the $\bp$-weighted version of $L_H$, the \emph{skin conformal Laplacian} $\bp^{-2} \cdot L_H$. \\

From $\tau(\bp,H)  \cdot \frac{k- 2}{4 (k-1)} \le  \frac{k- 2}{4 (k-1)}$, we notice:
{\small \begin{equation}\label{esti}
 \lambda^{\bp}_{H} \le  \frac{k- 2}{4 (k-1)} < 1/4.
\end{equation}}

In the regular case where $\Sigma \v$, $\lambda^{\bp}_{H}$ clearly is the usual first eigenvalue of $\bp^{-2} \cdot L_H$ and, different from the singular case we discuss below, there are no non-trivial solutions of the eigenvalue equation for any $\lambda < \lambda^{\bp}_{H}$.\\

\subsubsection{Regularizations and Subcritical Operators} \label{sus}
\bigskip

Since $\bp$ is locally Lipschitz regular we note for any convergent sequence of minimal hypersurfaces $H_i$ a $C^\alpha$-convergence of the $\bp_{H_i}$. Also, when $H$ is sufficiently regular, the eigenfunctions of $\bp^{-2} \cdot L_H$ are $C^{2,\alpha}$-smooth, for
$\alpha  \in (0,1)$.\\

Now, we observe, that, in the procedure of inductive descent, chosen in this paper, we lose some regularity during the dimensional descent. But this issue can easily be fixed. The distinctive geometry of our setup allows us to restore a better regularity by means of a simple smoothing of $\bp$, we repeatedly apply in each inductive loop.\\

To define this smoothing of $\bp$, we choose a cut-off function $\phi_\varrho(z)=\psi(|F_m(z)|/\varrho)$, for some $\psi \in C^{\infty}(\R,[0,1])$, with $\psi \equiv 0$ on $\R^{\ge 1}$ and $\psi \equiv 1$ on $\R^{\le 1/2}$, for $10 \cdot \gamma \le  \varrho$, for some $\varrho \in (0,1/10)$.\\

Now we modify $\bp_X$ on the torus component
\begin{equation}\label{xt}
XT_\varrho(\ve,\gamma):=X(\ve,\gamma) \cap QT^m_{\varrho},
\end{equation}
 so that it becomes constant on $XT$, being equal the value of $\bp$ in $F_m^{-1}(1_T)$:
\[\bsp_X(z) := \phi_\varrho \cdot \bp_X(z) + (1-\phi_\varrho) \cdot \bp_X(F_m^{-1}(1_T)),  \mm{ for } z  \in X  \setminus \Sigma_X.\]
Thus, on $XT_\varrho$,  keeping $\varrho>0$ fixed, $\bsp_X$ is smooth. Since $\bp_X(F_m^{-1}(1_T)) \ra  0$, for $\ve,\gamma \ra 0$, we trivially have:
\begin{equation}\label{bss}
 \bsp_X(\ve,\gamma) \ra  0, \mm{ for } \ve,\gamma \ra 0, \mm{ in } C^j\mm{-norm, on } XT_\varrho, \mm{ for any } j >0.
 \end{equation}

  To understand the eigenvalues of $\bsp_X^{-2} \cdot L_X$ we note the following estimate.

\begin{lemma}\label{lhh0}  \, For sufficiently small $\ve, \gamma>0$, we have   \[1/2 \cdot\bp_X(z) \le \bsp_X(z) \le 2 \cdot\bp_X(z), \mm{ for any } z  \in X  \setminus \Sigma_X.\]
Thus, for such $\ve, \gamma>0$, the generalized principal eigenvalue $\lambda^{\bsp}_{X}$ of $\bsp_X^{-2} \cdot L_X$ satisfies:
{\small \begin{equation}\label{esti2}
 0 < 1/4 \cdot \lambda^{\bp}_{X} \le \lambda^{\bsp}_{X} \le 4 \cdot \lambda^{\bp}_{X} < 1/2.
\end{equation}}
\end{lemma}

\textbf{Proof} \quad We recall the following form of a Harnack  inequality for $\bp$, one easily derives from the skin transform axioms, [L3],Lemma 3.2(ii):
\begin{equation}\label{liph} \left|\bp(x)/\bp(p)  - 1\right|  \le 2 \cdot L \cdot \bp(p) \cdot d(x,p)\end{equation}  for any $q \in B_{1/(2 \cdot L \cdot \bp(p))}(p)$. For  $\ve, \gamma>0$ small enough we know from \ref{rs} that $XT$ is arbitrarily  close to flat metric and thus $\bp$ arbitrarily close to $0$. Since the diameter of $XT$ is uniformly bounded, for $\ve,\gamma \ra 0$, we may assume $XT \subset B_{1/(2 \cdot L \cdot \bp(F_m^{-1}(1_T)))}(F_m^{-1}(1_T))$.\\

 For any $\eta >0$, we infer from (\ref{liph}), for  sufficiently small $\bp(F_m^{-1}(1_T))$,
that \[(1- \eta) \cdot \bp(F_m^{-1}(1_T))  \le \bp(z) \le (1+ \eta) \cdot \bp(F_m^{-1}(1_T)), \mm{ for any } z  \in XT.\]
This gives the claims for sufficiently small  $\ve, \gamma>0$. \qed\\

\subsubsection{Martin Integrals} \label{sus}
\bigskip

 Now we reach some of the essential analytic features of the operators $L_X$ and $\bsp_X^{-2} \cdot L_X$ on the non-compact space $X \setminus \Sigma$, we get in the case where $X$ is singular.\\

 Their spectral theory clearly differs substantially from that on closed smooth manifolds. In the singular case we have a criticality theory that suggests to define associated skin adapted operators $L_{X,\lambda}$, for any $\lambda <\lambda^{\bp}_{X}$.

\begin{proposition}\emph{\textbf{(Subcritical Variants of $L_X$)}}\label{lh1} \, For $m>3$, we consider the operators
 \[L_{X,[\lambda]}:= - \triangle   + \frac{m-3}{4 (m-2)}  \cdot scal_X  -\lambda  \cdot  \bsp_X^2, \mm{ for any } \lambda <\lambda^{\bsp}_{H}/8.\]
Now we assume that $\Sigma \n$. Then the  $L_{X,[\lambda]}$ are skin adapted and, from this  we get
 \begin{enumerate}
   \item The Martin boundary of $L_{X,[\lambda]}$ is homeomorphic to $\Sigma$. In particular, a function $u>0$ solves $L_{X,[\lambda]} \, v = 0$ if and only if $u$ admits a integral representation in terms of a (uniquely determined) finite Radon measure $\mu=\mu_u$ on $\Sigma$:
   \begin{equation}\label{yy}
       u(x) = \int_{\Sigma} k_{L_{X,[\lambda]}}(x;y) \, d \mu(y), \, \mm{ on } X \setminus \Sigma,
\end{equation}
 where $k_{L_{X,[\lambda]}}(x;y)$ denotes the Martin kernel of $L_{X,[\lambda]}$ on $X \setminus \Sigma$.
   \item For $\lambda \in (0,\lambda^{\bsp}_{X}/8)$, every solution of $L_{X,[\lambda]} \, w=0$, is $C^{2,\alpha}$-smooth on $X \setminus \Sigma$ and it is smooth on $XT$. Thus $scal(u^{4/m-3} \cdot g_X)$ is well-defined and we have
       \begin{equation}\label{yyy}
       scal(u^{4/m-3} \cdot g_X) \cdot u^{\frac{m+1}{m-3}} =  \gamma_{m-1} \cdot L_X \,   u =   \gamma_{m-1} \cdot\lambda  \cdot  \bsp_X^2  \cdot u >0.
       \end{equation}
 \end{enumerate}
\end{proposition}

\textbf{Proof} \quad The first skin adaptedness condition (i) follows from \ref{lh0} and \ref{lhh0}. For the second condition (ii) we have from \ref{lhh0} for any smooth function $f$ on $X$ with $supp \: f \subset H \setminus \Sigma$:
\begin{equation}\label{cf1}
\int_X  f  \cdot  L_X f  \,  dA  \ge \lambda^{\bp}_X\cdot \int_X \bp^2 \cdot f^2   \, d A \ge
\end{equation}
\[ \lambda^{\bp}_X/2\cdot \int_X \bp^2 \cdot f^2   \, d A  + \lambda^{\bsp}_X/8\cdot \int_X \bsp^2 \cdot f^2   \, d A\]
Therefore, condition (ii) is satisfied for $\lambda <\lambda^{\bsp}_X/8$:
\begin{equation}\label{cf2}
\int_X  f  \cdot  L_{X,[\lambda]} f  \,  dA  \ge \lambda^{\bp}_X/2\cdot \int_X \bp^2 \cdot f^2   \, d A.
\end{equation}
(Of course, \ref{lhh0} and this argument leave space to improve the estimates, but here we merely need this result for very small positive $\lambda$.)\\

From the skin adaptedness of $L_{X,[\lambda]}$, one may derive assertion (i). For this, we refer to [L4],Th.5 and Th.9 (and Ch.5.1.),  based on Th.10.\qed

In other words, for  $\Sigma \n$, we may  use conformal deformations by eigenfunctions for $\lambda>0$ strictly smaller than $\lambda^{\bsp}_X/8$ to get a $scal >0$-geometry on $X \setminus \Sigma$, while we can take profit from the asymptotic analysis, for instance, the Martin theory of (i), for skin adapted operators. In fact, in our further discussion, we will use a rather small $\lambda>0$, to ensure uniform estimates for the growth of the Martin kernel towards $\Sigma$.\\

\subsubsection{Horizons and Torus Components} \label{mos}
\bigskip

We define and analyze  conformal deformations of the area minimizing hypersurfaces $X^{m-1} \subset Q^m$, obtained in Ch.\ref{ami}, to get $scal >0$-metrics with the two additional properties that the deformations give a horizon around the singular set and they keep the torus component nearly flat.\\

 \textbf{Singular Case} \, When $\Sigma\n$, the idea is apply \ref{lh1} for an \emph{equidistributed} positive Radon measure $\mu$ on $\Sigma_X$ and to use growth estimates for the Martin kernel, we get for a small $\lambda >0$, to show that $u_\mu(x) = \int_{\Sigma} k_{L_{X,[\lambda]}}(x;y) \, d \mu(y)$
is a solution that induces an inner horizon shielding $\Sigma$.

\begin{remark} \textbf{(Gravitating Mass of $\Sigma$)} \, Although the existence of these $X$ will eventually be brought to a contradiction, we observe a detail valid independent of this particular
context.\\

When we think of $X$, after this deformation, as a space (partially) representing the original gravitational system, we notice an appealing interpretation: even the hypersurface singularities do not become visible as naked singularities, in the sense of Penrose [Pe], explained in a series of arguments in favor of a cosmic censorship. Also, in the spirit of the Penrose inequality, we observe, that hypersurface singularities give rise to positive mass contributions, which can be  lower estimated in terms of the area of inner  horizons shielding $\Sigma$.\qed
\end{remark}

In the actual technical implementation of these ideas, in [L4], we use an approximation scheme since the potential complexity $\Sigma$ makes it difficult to define and control such a $\mu$ directly. The outcome is a deformation function build in two steps:
 \begin{itemize}
   \item For any $\delta >0$, we choose a point $p_\delta \in X \setminus \Sigma$  with $dist(p_\delta,\Sigma)=\delta$ and take the minimal Green's function $G(\cdot,p_\delta)>0$ of $L_{X,[\lambda]}$ with pole in $p_\delta$ as a basic solution of minimal growth along $\Sigma$. For such solutions, the skin adaptedness of $L_{X,[\lambda]}$ grants a good control over the asymptotic behavior along $\Sigma$. We note that, also due to the skin adaptedness,
       there are no regular positive solutions on $X \setminus \Sigma$ of minimal growth along $\Sigma$.
   \item The bending effect of $G(\cdot,p_\delta)$ towards $\Sigma$ is not strong enough to give a repelling horizon. Therefore, we append a second deformation. In a $\delta^*$-neighborhood $U_{\delta^*}(\Sigma)$ of $\Sigma$, for $\delta \gg \delta^*>0$, we modify $G(\cdot,p_\delta)$ and merge $G(\cdot,p_\delta)$ with localized solutions associated to some piecewise equidistributed positive Radon measures along some approximation of $\Sigma$.
 \end{itemize}
We denote the resulting positive function by  $\psi[\delta,\delta^*]$.\\

   Since $L_{X,[\lambda]}=- \triangle   + \frac{m-3}{4 (m-2)}  \cdot scal_X  -\lambda  \cdot  \bsp_X^2$ is skin adapted and it coincides with  $- \triangle   + \frac{m-3}{4 (m-2)}  \cdot scal_X  -\lambda  \cdot  \bp_X^2$ near $\Sigma_X$,
 we can apply [L5], Th.1, explicitly stated in Prop.4.15 and 4.17, to build the following deformed metric on $X^{m-1}$:

\begin{proposition}\emph{\textbf{(Horizon Metrics)}}\label{hmw} \, For any  $\delta \gg \delta^*>0$ and sufficiently small $\lambda>0$, $\lambda \in (0,\lambda^{\bsp}_X/8)$,  we can define the horizon metric $g[\delta,\delta^*]:=\psi[\delta,\delta^*]^{4/(m-3)}\cdot g_{X\setminus \Sigma_X}$ with the following properties:
\begin{enumerate}
\item On $X^{m-1}\setminus U_{2 \cdot \delta}(\Sigma)$, the function $\psi[\delta,\delta^*]$ is a bounded solution of $L_{X,[\lambda]} \, w=0$:
\begin{equation}\label{psi}
- \triangle \,   \psi   + \frac{m-3}{4 (m-2)}  \cdot scal_X \cdot \psi= \lambda \cdot  \bsp_X^2  \cdot \psi
\end{equation}
\item For balls $B_{\delta^+}(p_\delta)$, for small $\delta^+>0$, $\delta^+ \ll |\delta-\delta^*|$, the distance sphere $\p B_{\delta^+}(p_\delta)$, relative
$g_{X\setminus \Sigma}$, is smooth and has positive mean curvature, relative $g[\delta,\delta^*]$.
\item There is a smoothly bounded neighborhood $V \subset U_{\delta^*}(\Sigma)$ of $\Sigma$, so that $\p V$ has positive mean
    curvature, relative $g[\delta,\delta^*]$, and $scal (g[\delta,\delta^*]) > 0$ on $X\setminus V$.
\end{enumerate}
\end{proposition}
For (ii) we note that, the conformal deformation of $g_{H\setminus \Sigma}$ with $G(\cdot,p_\delta)$ leads to an asymptotically flat end around  $p_\delta$. This causes the
 positive mean curvature of small distance spheres. Hence, area minimizers in the deformed geometry on $X\setminus V$ will not approach $\p B_{\delta^+}(p_\delta)$ but be repelled. \\

\textbf{Regular Case} \, When $X$ is a closed smooth hypersurface, we choose the first eigenfunction $\psi_1>0$ associated to $\lambda^{\bsp}_X>0$ of
${\bsp}^{-2}_X \cdot L_X$. It solves $L_{X,[\lambda^{\bsp}_{X}]} \, w=0$:
\begin{equation}\label{rrc}
- \triangle \,   \psi_1   +\frac{m-3}{4 (m-2)} \cdot scal_X  \cdot \psi_1= \lambda^{\bsp}_{X}  \cdot  \bsp_X^2  \cdot \psi_1
\end{equation}
From small $C^n$-perturbations of the ambient space, we can always assume that $X$ is not totally geodesic and therefore $\bp_X$ and, hence, $ \bsp_X$ are strictly positive.\\

Thus, in this case, we take the $scal>0$-metric $\psi_1^{4/m-3} \cdot g_X$.\\

\textbf{Recovery of Torus Components} \, Now we recall from the proof of \ref{rs}, that for any given $\varrho>0$, for $\ve \ra 0$ and $\gamma \ra 0$:  $XT_{\varrho/2}(\ve,\gamma)$ converge to $T^{m-1}[0_T] \setminus B_{\varrho/2}(0_T)$ in $C^{l^*}$-topology in the sense that they can be represented as graphs over $T^{m-1}[0_T] \setminus B_{\varrho/2}(0_T)$
and the associated functions converge to zero.\\

We henceforth also use these well-defined smooth parameterizations of the $XT_{\varrho/2}(\ve,\gamma)$ by $T^{m-1}[0_T] \setminus B_{\varrho/2}(0_T)$  to compare
functions on varying $XT_{\varrho/2}(\ve,\gamma)$, in [L3]-[L5] this canonical parametrization is called the identification map $\D$.\\

 With these parametrizations, the $C^{l^*}$-convergence of the $XT_{\varrho/2}(\ve,\gamma)$ to a flat limit gives us a $C^{l^*-2}$-convergence of the coefficients  of our operators, for $\ve \ra 0$, $\gamma \ra 0$:
\[scal_{X} \circ F_m^{-1} \ra 0, \,  |A_{X} | \circ F_m^{-1}  \ra 0 \mm{ and } \bsp_{X}(\ve,\gamma) \circ F_m^{-1} \ra 0,\]
on $T^{m-1}[0_T] \setminus B_{\varrho/2}$.\\

Moreover, the estimate (\ref{esti}) shows that in both cases $\lambda$ respectively $\lambda^{\bsp}_X$ belong to the bounded interval $(0,2)$.
Therefore, we have common Harnack inequalities and elliptic regularity estimates with uniform estimates, for $\ve \ra 0$ and $\gamma \ra 0$.\\

We infer that, for any sequence of $\ve_k \ra 0,\gamma_k \ra 0$,  for $k \ra \infty$, and any sequence of solutions  $\psi(\ve_k,\gamma_k)>0$, which may either belong to the singular case  (\ref{psi}) or to the regular case (\ref{rrc}), defined on $XT_{\varrho/2}(\ve_k,\gamma_k)$, normalized to $\psi(F_m^{-1}(1_T))=1$,   $C^{l^*-1}$-uniform subconvergences to a positive harmonic function $u_\varrho$  on $T^{m-1}[0_T] \setminus B_{3 \cdot \varrho/4}(0_T)$.\\

 (The regularity theory for $C^{k,\alpha}$-regular operators and solutions would provide  better  estimates, but we do not need them in this discussion.)\\

Now we let also $\varrho \ra 0$, then a diagonal sequence argument shows that the limit is a positive harmonic $u$  on $T^{m-1} \setminus \{0_T\}$. However, applying the Harnack inequality and the maximum principle one easily shows that every positive harmonic $u$  on $T^{m-1} \setminus \{0_T\}$ is constant. Summarizing, we get

\begin{corollary}\emph{\textbf{(Torus Components)}}\label{hm2} \, For $m>3$ and any $\varrho> 0$, $\xi >0$, there are some small $\ve > 0,\gamma > 0$, so that
 $(XT_\varrho(\ve,\gamma),\psi^{4/m-3} \cdot g_X)$ is  almost isometric to the flat $T^{m-1}[0_T] \setminus B_\varrho(0_T)$:
\[|\psi^{4/m-3} \cdot g_X-F^*_m(g_{flat})|_{C^{l^*-1}(XT_\varrho(\ve,\gamma))} < \xi.\]
\end{corollary}
\medskip

\textbf{Conclusion} \, Of course, since we are free to choose $l^*$ as large as we want, this allows us to complete the \textbf{proof} of the main inductive descent result \ref{id} until reach the three-dimensional case: we combine  \ref{hmw}, in case $X$ is singular, respectively the choice of (\ref{rrc}), when $X$ is regular, with the latter result \ref{hm2}. \\

Now, we reach the final step in the tower of inductive descent where we pass from $m=3$ to a surface. As explained in \ref{teo}, for the proof of Theorems, we merely need to see that the resulting closed surface $X^2$ has genus $\ge 1$, this follows from \ref{rs}, and the integral scalar curvature of $X^2$ is positive. To this end we notice, from the Gau\ss-Codazzi equation, the
stability inequality $Area''(f) \ge 0$ in (\ref{stabi})  can now be written
\begin{equation}\label{stl2}
 \int_{X}|\nabla_X f|^2 - \left(1/2 \cdot Scal_{Q^3}  + 1/2 \cdot |A|^2 - Scal_X \right) \cdot f^2 \: dVol \ge 0
\end{equation}
for any smooth function $f$ on $X$. We choose $f \equiv 1$ to see that
\begin{equation}\label{stl2}
 \int_{X} Scal_X \: dVol \ge  1/2 \cdot \int_{X} Scal_{Q^3}  +  |A|^2 \: dVol > 0
\end{equation}

(We leave it to the interested reader to formally derive the other conditions we formulated in \ref{id} also in the surface case, using the  two dimensional transformation law $scal(e^{2 \cdot u} \cdot g_{X^2}) \cdot e^{2 \cdot u} =-\Delta\,  u  + scal(g_{X^2})$.) \qed

\subsubsection{Vanishing Energy Case} \label{bo}
\bigskip

Finally, we consider the rigidity assertion of Theorem 1: $E(M,g) = 0$ if and only if $(M,g)$ is isometric to the Euclidean space $\R^n$. To derive this result in arbitrary dimensions we apply the main part of Theorem 1 we have already proved above: for any asymptotically flat $(M^n,g)$ of order $\tau  > \frac{n-2}{2}$ with $scal_M\ge 0$ we have $E(M,g) \ge 0$.\\

With this non-negativity of $E$, actually applied two times, the classical subsequent arguments apply regardless of the dimension. For the sake of completeness, we sketch a version of these arguments and refer to [LP],10.7 for details.\\

\textbf{Outline} \, If $E(M,g) = 0$, then $(M,g)$ must be scalar flat. Otherwise we could find a conformal deformation to another such manifold, but with $E<0$, reversing the argument we used to reduce Theorem 1 to the non-existence of a $scal>0$-island of [L1],Ch.6. This contradicts the proven non-negativity of $E$. \\

 Now we show that such an $(M,g)$ must be Ricci flat. The trick is to choose, around any point $p \in M$, a perturbation $\phi \cdot Ric_M$, where $\phi \ge 0$ is a smooth compactly supported function with $\phi(p)=1$.
 Since this perturbation is compactly supported, we still have $E(M,g+t \cdot \phi \cdot Ric_M)=0$, for any $t \in \R$.\\

 One can show that for small $|t|>0$,  $g+t \cdot \phi \cdot Ric_M$ can be conformally deformed into another scalar-flat and asymptotically flat manifold $(M,g_t)$. Again we know that $E(M,g_t)\ge 0$. We use the variational formula for the energy, for $scal(M,g)\equiv 0$:
 \[\frac{d}{dt}E(M,g_t)(0)=(n-2)^{-1} \cdot \int_M \langle Ric_M, \phi \cdot Ric_M\rangle \, dV.\]
Since $E(M,g_0)= 0$ the real function $E(M,g_t)$ attains its minimum for $t=0$,  this derivative vanishes and we infer that around any given point $Ric_M \equiv 0$.\\

Thus, when $E(M,g) = 0$, then $(M,g)$ must even be Ricci flat. Next, the asymptotically flatness expressed in harmonic coordinates allows us to exploit the Ricci flatness by a Bochner argument. The outcome are $n$ orthogonal globally defined parallel vector fields yielding another upgrade, now from Ricci flatness, to proper flatness.\qed


{\small
\begin{thebibliography}{LLLL}
\bibitem[ADM]{ADM} Arnowitt, R., Deser, S. and Misner, C.: Coordinate invariance and energy expressions in general relativity, Phys.Rev.122 (1961), 997-1006
\bibitem[B]{B} Bartnik, R.: The Mass of an Asymptotically Flat Manifold, Comm. Pure Appl. Math. 34 (1986), 661-693
\bibitem[E]{E} Eichmair, M.: The Jang equation reduction of the spacetime positive energy theorem in
dimensions less than eight. Comm. Math. Phys. 319 (2013), 575-593
\bibitem[EHLS]{EHLS} Eichmair, M., Huang, L.-H. Lee, D. and Schoen, R.: The spacetime positive mass theorem in dimensions less than eight, J. Eur. Math. Soc. 18 (2016), 83-121
\bibitem[F]{F} Federer, H.: Geometric Measure Theory, Spinger Verlag, Berlin  (1969)
\bibitem[GMS]{GMS} Giaquinta, M., Modica, G. and Sourek, J.: Cartesian Currents in the Calculus of Variations, Vol. I, Springer (1998)
\bibitem[KMS]{KMS} Khuri, M., Marques, F. and Schoen, R.: A compactness theorem for the Yamabe problem, J. Diff. Geom. 81 (2009), 143-196
\bibitem[KP]{KP} Krantz, S. and Parks, H.: Geometric Integration Theory, Birkh\"auser, Boston (2008)
\bibitem[LP]{LP} Lee, J. and T. Parker, T.: The Yamabe problem, Bull. AMS 17 (1987), 37-91
\bibitem[L1]{L1} Lohkamp, J.: Scalar Curvature and Hammocks, Math. Ann. 313 (1999), 385-407
\bibitem[L2]{L2} Lohkamp, J.: Curvature h-principles, Ann. of Math. 142 (1995), 457-498
\bibitem[L3]{L3} Lohkamp, J.: Skin Structures on Minimal Hypersurfaces (2015), arXiv:1512.08249
\bibitem[L4]{L4} Lohkamp, J.: Hyperbolic Geometry and Potential Theory on Minimal Hypersurfaces (2015), arXiv:1512.08251
\bibitem[L5]{L5} Lohkamp, J.:  Skin Structures in Scalar Curvature Geometry (2015),  arXiv: 1512.08252
\bibitem[L6]{L6} Lohkamp, J.: The Higher Dimensional Positive Mass Theorem II, Preprint
\bibitem[P]{P} Pitts, J.: Existence and Regularity of Minimal Surfaces on Riemannian Manifolds, Math. Notes, Princeton Univ. Press (1981)
\bibitem[Pe]{Pe} Penrose, R.: Naked Singularities, Ann. N.Y . Acad. Sci. 224 (1973), 125-134.
\bibitem[PT]{PT} Parker T., Taubes C.: On Witten's proof of the
positive energy theorem, Comm. Math. Phys. 84 (1982), 223-238
\bibitem[S]{S} Schoen, R.: Variational theory for the Total Scalar curvature Functional for Riemannian Metrics and Related Topics , in Topics in Calculus of Variations,
LNM 1365 , Springer  (1989),120-154
\bibitem[SY1]{SY1} Schoen, R. and Yau, S.T.: Existence of
incompressible minimal surfaces and the topology of three dimensional manifolds with non-negative scalar curvature, Ann. of Math. 110 (1979), 127-142
\bibitem[SY2]{SY2} Schoen, R. and Yau, S.T.: On the proof of the positive mass conjecture in general relativity, Comm. Math. Phys. 65 (1979), 45-76
\bibitem[SY3]{SY3} Schoen, R. and Yau, S.T.: Proof of the positive mass theorem II, Commun.Math.Phys.79 (1981),231-260
\bibitem[Si1]{Si1} Simon, L.: Lectures on Geometric Measure Theory,
Proceedings of the Centre for Mathematical Analysis, Australian National University, Canberra, 1983
\bibitem[Si2]{Si2} Simon, L.: A Strict Maximum Principle for Area-Minimizing Hypersurfaces, J. Diff. Geom.
26 (1987), 327-335
\bibitem[W]{W} Witten, E.: A new proof of the positive energy theorem, Comm. Math. Phys. 80 (1981),
381-402
\end{thebibliography}}
\end{document}